%
%
\input amstex.tex
\documentstyle{amsppt}
\magnification=1200
\baselineskip=13pt
\hsize=6.5truein
\vsize=8.9truein
\countdef\sectionno=1
\countdef\eqnumber=10
\countdef\theoremno=11
\countdef\countrefno=12
\countdef\cntsubsecno=13
\sectionno=0
\def\newsection{\global\advance\sectionno by 1
                \global\eqnumber=1
                \global\theoremno=1
                \global\cntsubsecno=0
                \the\sectionno}

\def\newsubsection#1{\global\advance\cntsubsecno by 1
                     \xdef#1{{\S\the\sectionno.\the\cntsubsecno}}
                     \ \the\sectionno.\the\cntsubsecno.}

\def\theoremname#1{\the\sectionno.\the\theoremno
                   \xdef#1{{\the\sectionno.\the\theoremno}}
                   \global\advance\theoremno by 1}

\def\eqname#1{\the\sectionno.\the\eqnumber
              \xdef#1{{\the\sectionno.\the\eqnumber}}
              \global\advance\eqnumber by 1}

\def\thmref#1{#1}

\global\countrefno=1

\def\refno#1{\xdef#1{{\the\countrefno}}\global\advance\countrefno by 1}


\def\C{{\Bbb C}}
\def\Z{{\Bbb Z}}
\def\Zp{{\Bbb Z}_+}
\def\hf{{1\over 2}}

\def\A{{\Cal A}_q}
\def\a{\alpha}
\def\b{\beta}

\def\d{\delta}

\def\th{\theta}

\def\m{\mu}
\def\n{\nu}

\def\vp{\varphi}

\refno{\AlSaC}
\refno{\AskeI}
\refno{\AskeW}
\refno{\FlorVeen}
\refno{\FlorVtwee}
\refno{\GaspR}
\refno{\IsmaW}
\refno{\KalnMaMi}
\refno{\KalnM}
\refno{\KalnMMeen}
\refno{\KalnMMtwee}
\refno{\KoelJMAA}
\refno{\KoelITSF}
\refno{\KoelDMJ}
\refno{\KoelIM}
\refno{\KoelCJM}
\refno{\KoelS}
\refno{\KoorS}
\refno{\Moak}
\refno{\Rahm}
\refno{\Swar}
\refno{\VaksK}
\refno{\VanAK}
\refno{\VileK}
\refno{\Wats}
\topmatter
\title Yet another basic analogue of Graf's addition formula \endtitle
\rightheadtext{Basic analogue of Graf's addition formula}
\author H.T. Koelink\endauthor
\affil  Katholieke Universiteit Leuven\endaffil
\address Departement Wiskunde, Katholieke Universiteit Leuven,
Celestijnenlaan 200 B, B-3001 Leuven (Heverlee), Belgium\endaddress
\email erik\%twi\%wis\@cc3.KULeuven.ac.be\endemail
\date September 20, 1994\enddate
\thanks Supported by a Fellowship of the Research Council of
the Katholieke Universiteit Leuven. \endthanks
\abstract
An identity involving basic Bessel functions and Al-Salam--Chihara
polynomials is proved for which we recover Graf's addition formula for
the Bessel function  as the base $q$ tends to $1$. The corresponding product
formula is derived. Some known
identities for Jackson's $q$-Bessel functions are obtained as limiting cases.
As special cases we prove identities for $q$-Charlier polynomials.
\endabstract
\endtopmatter
\document


\subhead\newsection . Introduction and formulation of results\endsubhead
A classical result for the Bessel function $J_\n(z)$ of order $\n$ and
argument $z$ defined by the absolutely convergent series
$$
J_\nu(z) = \sum_{k=0}^\infty {{(-1)^k(z/2)^{\nu+2k}}\over{k!\,
\Gamma(\nu+k+1)}}
$$
is the addition formula
$$
J_\nu\Bigl( \sqrt{x^2+y^2-2xy\cos \psi}\Bigr)
\left( {{x-ye^{-i\psi}}\over{x-ye^{i\psi}}} \right)^{\nu/2} =
\sum_{m=-\infty}^\infty J_{\nu+m}(x)J_m(y)e^{im\psi},
\tag\eqname{\vglGrafadditie}
$$
$\vert ye^{\pm i\psi} \vert < \vert x\vert$,
due to Graf (1893), cf. \cite{\Wats , \S 11.3(1)}, for general $\nu$, and
due to Neumann (1867) for $\nu=0$, cf. \cite{\Wats , \S 11.2(1)}.
In case $\nu\in\Z$ the conditions on $x$, $y$ in \thetag{\vglGrafadditie}
can be removed.
The corresponding product formula for the Bessel function is
$$
J_{\nu+m}(x)J_m(y)=
{{1}\over{2\pi}} \int_0^{2\pi}
J_\nu\Bigl( \sqrt{x^2+y^2-2xy\cos \psi}\Bigr)
\left( {{x-ye^{-i\psi}}\over{x-ye^{i\psi}}} \right)^{\nu/2}
e^{-im\psi}\,d\psi .
\tag\eqname{\vglGrafproduct}
$$

There are many $q$-analogues of the Bessel function and for several of these
$q$-analogues of the Bessel function there exist identities which have the
Graf addition formula \thetag{\vglGrafadditie} as the limit for $q\uparrow 1$.
The structure of these identities may be very different from the structure of
\thetag{\vglGrafadditie}. These $q$-analogues of the Graf addition formula
often follow from a certain interpretation of a $q$-Bessel function on
quantum groups \cite{\KoelDMJ}, \cite{\KoelIM}, \cite{\VaksK}, or on quantum
algebras \cite{\FlorVeen}, \cite{\KalnM}, \cite{\KalnMMtwee}, or from
a certain generating function for a $q$-Bessel function \cite{\KoorS},
\cite{\KoelITSF}, a method closely related to the quantum algebra approach.
The first two methods are motivated by the group theoretic
proof of the Graf addition formula as given by Vilenkin and Klimyk
\cite{\VileK, \S 4.1.4(2)}.

It is the purpose of this note to give analytic proofs of a $q$-analogue
of Graf's addition formula, which was first obtained by a formal calculation
using the quantum group of plane motions, and of the corresponding product
formula. In the rest of this introduction we formulate the addition and
product formula. In the next section we prove the
product formula using a connection
coefficient formula for the Al-Salam--Chihara polynomials. In \S 3 we derive
the addition formula. In \S 4 we consider some special and limiting cases,
and in particular the limit case $q\uparrow 1$,
and we present some links with known results in this direction. In this
section we also derive identities for $q$-Charlier
polynomials as special cases. Finally, in \S 5 we
remark very shortly on the link with the quantum group of plane motions.

In order to formulate the results we recall the notation for the
$q$-shifted factorial;
$$
(a;q)_k=\prod_{i=0}^{k-1}(1-aq^k),\quad (a_1,\ldots,a_r;q)_k =
(a_1;q)_k\ldots (a_r;q)_k,
\quad k\in\Zp\cup\{\infty\}
$$
and for the $q$-hypergeometric series
$$
\align
{}_r\vp_s\left({{a_1,\ldots,a_r}\atop{b_1,\ldots,b_s}};q,z\right) &=
{}_r\vp_s(a_1,\ldots,a_r;b_1,\ldots,b_s;q,z) \\
&= \sum_{k=0}^\infty {{(a_1,\ldots,a_r;q)_k}\over{(q,b_1,\ldots,b_s;q)_k}}
\Bigl( (-1)^k q^{\hf k(k-1)}\Bigr)^{1-r+s} z^k.
\endalign
$$
We always assume $0<q<1$. These notations follow the book \cite{\GaspR}
by Gasper and Rahman, which should be consulted for more information on
this subject. For generic values of the parameters the region of
convergence of this series is $\infty$, $1$, or $0$, according to
$r<s+1$, $r=s+1$, or $r>s+1$.

Note that a $q$-hypergeometric series with $q^{-n}$, $n\in\Zp$, as one
of the upper parameters terminates, since $(q^{-n};q)_k=0$ for
$k>n$. If $q^{-n}$ occurs as a lower parameter, then the $q$-hypergeometric
series is in general not well-defined.
But for $n\in\Zp$ we use the following convention
$$
(q^{1-n};q)_\infty \sum_{k=0}^\infty {{c_k}\over{(q^{1-n},q;q)_k}} =
\sum_{k=n}^\infty c_k {{(q^{1-n+k};q)_\infty}\over{(q;q)_k}}
= (q^{n+1};q)_\infty \sum_{k=0}^\infty {{c_{k+n}}\over{(q^{1+n},q;q)_k}}.
\tag\eqname{\vglnegnseries}
$$
(This corrects the first identity of the remark
following proposition~4.1 in \cite{\KoelIM}.)

In order to formulate the main result of this paper we introduce the
Al-Salam--Chihara polynomials;
$$
S_n(\cos\th;a,b\mid q) = a^{-n}(ab;q)_n \, {}_3\vp_2 \left(
{{q^{-n},ae^{i\th},ae^{-i\th}}\atop{ab,\; 0}};q,q\right).
\tag\eqname{\vgldefASCpols}
$$
These polynomials were originally introduced by Al-Salam and Chihara
\cite{\AlSaC} as the most general set of orthogonal polynomials
satisfying a certain `convolution property'. The orthogonality measure
for these polynomials has been obtained by Askey and Ismail
\cite{\AskeI, \S 3}, which is a special case of the more general
four-parameter class
of Askey-Wilson polynomials, cf. \cite{\AskeW}. The definition
\thetag{\vgldefASCpols} used here gives the Al-Salam--Chihara
polynomials as Askey-Wilson polynomials with two parameters set to zero.
The orthogonality relations for the
Al-Salam--Chihara polynomials are given by,
cf. Askey and Ismail \cite{\AskeI, \S 3.8},
Askey and Wilson \cite{\AskeW, thm.~2.2},
$$
{1\over{2\pi}} \int_0^\pi \bigl(S_kS_l\bigr)(\cos\th;a,b\mid q)
{{(e^{2i\th},e^{-2i\th};q)_\infty}\over{
(ae^{i\th},ae^{-i\th},be^{i\th},be^{-i\th};q)_\infty}}\, d\th =
{{\d_{k,l}}\over{(q^{k+1},abq^k;q)_\infty}},
\tag\eqname{\vglorthoASCpols}
$$
assuming that $\vert a\vert<1$, $\vert b\vert <1$. The Al-Salam--Chihara
polynomials are symmetric in the parameters $a$ and $b$, cf.
\cite{\AskeW, p.~6}.

The main result of this paper is the following addition formula, valid for
$\vert z\vert <1$, $\vert a\vert<1$, $\vert b\vert<1$;
$$
\multline
{{(q^{\n+1};q)_\infty}\over{(q;q)_\infty}}
{}_2\vp_1\left( {{ae^{i\th},ae^{-i\th}}\atop{q^{\n+1}}};q,z\right)
S_m(\cos\th;aq^{-\n},b\mid q) = \\
\sum_{n=-m}^\infty
(-1)^na^nz^nq^{\hf n(n-1)}
{{(q^{1+n};q)_\infty}\over{(q;q)_\infty}} \,
{}_1\vp_1\left({{q^{-m}}\atop{q^{1+n}}};q,a^2q^{m+n-\nu}z\right) \\
\times {{(q^{\n+n+1};q)_\infty}\over{(q;q)_\infty}} \,
{}_2\vp_1\left({{abq^{n+m},0}\atop{q^{\n+n+1}}};q,z\right)
S_{n+m}(\cos\th;a,b\mid q).
\endmultline
\tag\eqname{\vgladditionform}
$$
Of course, this formula can also be considered as a linearisation formula
for the product of a ${}_2\vp_1$-series and an Al-Salam--Chihara polynomial
in terms of another set of Al-Salam--Chihara polynomials.

The corresponding product formula is
$$
\multline
{{(q^{\n+1};q)_\infty}\over{2\pi (q;q)_\infty}}  \int_0^\pi
{}_2\vp_1\left( {{ae^{i\th},ae^{-i\th}}\atop{q^{\n+1}}};q,z\right)
S_m(\cos\th;aq^{-\n},b\mid q)\\  \times S_{n+m}(\cos\th;a,b\mid q)
{{(e^{2i\th},e^{-2i\th};q)_\infty}\over{
(ae^{i\th},ae^{-i\th},be^{i\th},be^{-i\th};q)_\infty}}\, d\th \\
= {{(-1)^na^nz^nq^{\hf n(n-1)}}\over{(q^{n+m+1},abq^{n+m};q)_\infty}}
{{(q^{1+n};q)_\infty}\over{(q;q)_\infty}} \,
{}_1\vp_1\left({{q^{-m}}\atop{q^{1+n}}};q,a^2q^{m+n-\nu}z\right) \\
\times {{(q^{\n+n+1};q)_\infty}\over{(q;q)_\infty}} \,
{}_2\vp_1\left({{abq^{n+m},0}\atop{q^{\n+n+1}}};q,z\right)
\endmultline
\tag\eqname{\vglproductform}
$$
for $\vert z\vert<1$, $\vert a\vert<1$, and $\vert b\vert <1$.
In this product formula we assume $\vert a\vert<1$, $\vert b\vert<1$,
but a similar product formula remains true for any choice of $a$ and $b$
for which the Al-Salam--Chihara polynomials are orthogonal polynomials.
In the general case only a finite number of discrete mass points have
to be added, cf. Askey and Wilson \cite{\AskeW, \S 2}.


\subhead\newsection . Proof of the product formula\endsubhead
In order to prove the product formula \thetag{\vglproductform} we
first prove two lemmas.

\proclaim{Lemma \theoremname{\lemone}} For the Al-Salam--Chihara
polynomials $S_n(\cdot;a,b\mid q)$ with $\vert a\vert<1$, $\vert
b\vert<1$ defined by \thetag{\vgldefASCpols}, $m,r\in\Zp$, $n\in\Z$ with
$n\geq -m$, we have
$$
\multline
{1\over{2\pi}}
\int_0^\pi S_m(\cos\th;aq^{-\n},b\mid q) S_{n+m}(\cos\th;a,b\mid q)
{{(ae^{i\th},ae^{-i\th};q)_r(e^{2i\th},e^{-2i\th};q)_\infty}\over{
(ae^{i\th},ae^{-i\th},be^{i\th},be^{-i\th};q)_\infty}}\, d\th \\
= (-a)^{-n} q^{n(\n+1)} q^{\hf n(n-1)}
{{(q^{\n+n+r+1};q)_\infty}\over{(q^{m+1},q^{\n+r+1},abq^{n+m+r};q)_\infty}}\\
\times {{(q^{1-n};q)_\infty}\over{(q;q)_\infty}}
\, {}_3\vp_2 \left( {{q^{-r},q^{-m-n},q^{-\n-n-r}}\atop{
q^{1-n},q^{1-m-n-r}/(ab)}};q,q^{r+1}{a\over b}\right).
\endmultline
$$
\endproclaim

\demo{Proof} We need the
connection coefficients for two sets of Al-Salam--Chihara polynomials
with one different parameter;
$$
S_n(x;\a,b\mid q) = \sum_{k=0}^n c_{k,n}(\a;a) S_k(x;a,b\mid q)
\tag\eqname{\vglccfASCpols}
$$
with
$$
c_{k,n}(\a;a)= {{(q^{-n};q)_k}\over{(q;q)_k}} a^{n-k}
(-1)^k q^{nk-\hf k(k-1)}(\a/a;q)_{n-k},
$$
given by Askey and Wilson \cite{\AskeW, (6.4), (6.5) with $c=d=0$}.

Now we start with the left hand side of the statement in the lemma. The
term $(ae^{i\th},ae^{-i\th};q)_r$ cancels part of the denominator of the
weight function. Next use \thetag{\vglccfASCpols} to write both
Al-Salam--Chihara polynomials in terms of Al-Salam--Chihara polynomials
with parameters $aq^r$ and $b$. Then we can use
\thetag{\vglorthoASCpols} to see that the left hand side of the
statement in the lemma equals
$$
\sum_{k=0}^m \sum_{l=0}^{n+m} \d_{k,l} {{c_{k,m}(aq^{-\n};aq^r)
c_{l,n+m}(a;aq^r)}\over{(q^{k+1},abq^{k+r};q)_\infty}}.
$$
If we use $k$ as the summation parameter, we can rewrite this as
$$
{{(q^{-\n-r};q)_m(q^{-r};q)_{n+m}}\over{(q,abq^r;q)_\infty}}
(aq^r)^{n+2m} \, {}_3\vp_2 \left( {{q^{-m},q^{-n-m},abq^r}\atop{
q^{\n+r+1-m},q^{1+r-n-m}}};q,q^{\n+2}a^{-2}\right).
$$
View this as a terminating ${}_3\vp_2$-series of degree $n+m$ to which
we apply the series inversion, cf. \cite{\GaspR, ex.~1.4(ii)},
$$
{}_3\vp_2\left( {{q^{-p},a,b}\atop{c,d}};q,z\right) =
{{(a,b;q)_p}\over{(b,c;q)_p}}(-z)^p q^{-\hf p(p+1)}
\, {}_3\vp_2\left({{q^{-p},q^{1-p}/c,q^{1-p}/d}\atop
{q^{1-p}/a,q^{1-p}/b}};q,{{cdq^{p+1}}\over{abz}}\right),
$$
$p\in\Zp$, to obtain a ${}_3\vp_2$-series as in the lemma. Some manipulations
with $q$-shifted factorials finish the proof of the lemma.
\qed\enddemo

The orthogonality relations \thetag{\vglorthoASCpols} for the
Al-Salam--Chihara polynomials show that the integral in lemma
\thmref{\lemone} is zero for $n>r$.
Observe that the right hand side of lemma \thmref{\lemone} also equals zero
for $n>r$ as follows from \thetag{\vglnegnseries} with $c_k=0$ for $k>n$
in this case.

The following lemma is straightforward and its proof is left to the
reader.

\proclaim{Lemma \theoremname{\lemtwo}} For $\vert dz\vert <1$ we have
$$
\multline
{{(q^{\m+1};q)_\infty}\over{(q;q)_\infty}}
\, {}_1\vp_1 \left( {a\atop{q^{\m+1}}};q,bz\right)
{{(q^{\n+1};q)_\infty}\over{(q;q)_\infty}} \, {}_2\vp_1
\left( {{c,0}\atop{q^{\n+1}}};q,dz\right)
= \\ {{(q^{\n+1};q)_\infty}\over{(q;q)_\infty}} \sum_{p=0}^\infty {{(dz)^p
(c;q)_p}\over{(q,q^{\n+1};q)_p}} {{(q^{\m+1};q)_\infty}\over{(q;q)_\infty}}
\, {}_3\vp_2 \left( {{q^{-p},q^{-p-\n},a}\atop{q^{\m+1},q^{1-p}/c}};q,
{{bq^{\n+p+1}}\over{dc}}\right)
\endmultline
$$
the last series being absolutely convergent.
\endproclaim

If we take $a=0$ in lemma \thmref{\lemtwo} and we replace $d$ by $d/c$ and
we let $c\to\infty$, then we essentially obtain the product formula for the
Hahn-Exton $q$-Bessel function, cf. Swarttouw \cite{\Swar, (3.1)}. If we
take $c=0$ in lemma \thmref{\lemtwo} and we replace $b$ by $b/a$ before
taking $a\to\infty$, then we obtain the product formula for Jackson's
$q$-Bessel functions, cf. Rahman \cite{\Rahm, (2.1)}.

The proof of the product formula \thetag{\vglproductform} can now be given.
Use the series representation for the ${}_2\vp_1$-series to see that
for $\vert z\vert<1$
$$
\multline
{{(q^{\n+1};q)_\infty}\over{2\pi (q;q)_\infty}}  \int_0^\pi
{}_2\vp_1\left( {{ae^{i\th},ae^{-i\th}}\atop{q^{\n+1}}};q,z\right)
S_m(\cos\th;aq^{-\n},b\mid q)\\  \times S_{n+m}(\cos\th;a,b\mid q)
{{(e^{2i\th},e^{-2i\th};q)_\infty}\over{
(ae^{i\th},ae^{-i\th},be^{i\th},be^{-i\th};q)_\infty}}\, d\th \\
={{(-a)^{-n} q^{n(\n+1)} q^{\hf n(n-1)}}\over{(q^{m+1},abq^{n+m};q)_\infty}}
{{(q^{\n+n+1};q)_\infty}\over{(q;q)_\infty}}  \sum_{r=0}^\infty
{{z^r(abq^{n+m};q)_r}\over{(q,q^{\n+n+1};q)_r}} \\
\times {{(q^{1-n};q)_\infty}\over{(q;q)_\infty}}
\, {}_3\vp_2 \left( {{q^{-r},q^{-m-n},q^{-\n-n-r}}\atop{
q^{1-n},q^{1-m-n-r}/(ab)}};q,q^{r+1}{a\over b}\right)
\endmultline
$$
by lemma \thmref{\lemone} and some straightforward simplifications.
Interchanging summation and integration is allowed for $\vert z\vert < 1$
by use of the estimate $\vert (ae^{i\th},ae^{-i\th};q)_r\vert \leq
(-\vert a\vert;q)_\infty^2$.
By lemma \thmref{\lemtwo} this expression is equal to
$$
\multline
{{(-a)^{-n} q^{n(\n+1)} q^{\hf n(n-1)}}\over{(q^{m+1},abq^{n+m};q)_\infty}}
{{(q^{1-n};q)_\infty}\over{(q;q)_\infty}} \,
{}_1\vp_1\left({{q^{-n-m}}\atop{q^{1-n}}};q,a^2q^{m-\nu}z\right) \\
\times {{(q^{\n+n+1};q)_\infty}\over{(q;q)_\infty}} \,
{}_2\vp_1\left({{abq^{n+m},0}\atop{q^{\n+n+1}}};q,z\right)
\endmultline
$$

Now use that \thetag{\vglnegnseries} implies that
$$
{{(q^{1-n};q)_\infty}\over{(q;q)_\infty}} \,
{}_1\vp_1\left({{a}\atop{q^{1-n}}};q,z\right) =
z^n (-1)^n q^{\hf n(n-1)}
{{(a,q^{1+n};q)_\infty}\over{(aq^n,q;q)_\infty}} \,
{}_1\vp_1\left({{aq^n}\atop{q^{1+n}}};q,zq^n\right),
$$
which holds for $n\in\Z$. This finishes the proof of the product formula
\thetag{\vglproductform}.


\subhead\newsection . Proof of the addition formula\endsubhead
The hard work for the proof of the addition formula
\thetag{\vgladditionform} has been done in the previous section. Let
$L^2\bigl( (-1,1),\, w(x)dx\bigr)$ denote the space of quadratically
integrable functions on $(-1,1)$ with respect to the weight function $w(x)$
defined by
$$
w(\cos\th) = {1\over{2\pi\sin\th}}
{{(e^{2i\th},e^{-2i\th};q)_\infty}\over{
(ae^{i\th},ae^{-i\th},be^{i\th},be^{-i\th};q)_\infty}}.
$$
The Al-Salam--Chihara polynomials $S_n(x;a,b\mid q)$ form a basis for this
$L^2$-space.

From the estimate $\vert(ae^{i\th},ae^{-i\th};q)_r\vert\leq (-\vert
a\vert;q)_\infty^2$ we obtain that the left hand side of
\thetag{\vgladditionform} as a function of $x=\cos\th$ is an element of this
$L^2$-space. Consequently, we may develop it in the basis of Al-Salam--Chihara
polynomials
$$
{{(q^{\n+1};q)_\infty}\over{(q;q)_\infty}}
{}_2\vp_1\left( {{ae^{i\th},ae^{-i\th}}\atop{q^{\n+1}}};q,z\right)
S_m(\cos\th;aq^{-\n},b\mid q) =
\sum_{n=-m}^\infty A_n
S_{n+m}(\cos\th;a,b\mid q)
$$
with
$$
\multline
A_n \int_{-1}^1
\bigl( S_{n+m}(x;a,b\mid q)\bigr)^2 w(x) \, dx =
{{(q^{\n+1};q)_\infty}\over{2\pi (q;q)_\infty}}  \int_0^\pi
{}_2\vp_1\left( {{ae^{i\th},ae^{-i\th}}\atop{q^{\n+1}}};q,z\right)
\\\times S_m(\cos\th;aq^{-\n},b\mid q) S_{n+m}(\cos\th;a,b\mid q)
{{(e^{2i\th},e^{-2i\th};q)_\infty}\over{
(ae^{i\th},ae^{-i\th},be^{i\th},be^{-i\th};q)_\infty}}\, d\th.
\endmultline
$$
Now use the orthogonality relations for the Al-Salam--Chihara polynomials
\thetag{\vglorthoASCpols} and the product formula
\thetag{\vglproductform} to find the correct value for $A_n$ as in the
addition formula \thetag{\vgladditionform}. This proves the
addition formula \thetag{\vgladditionform} as an identity in
$L^2\bigl( (-1,1),\, w(x)dx\bigr)$.

The left hand side of \thetag{\vgladditionform} is a continuous function of
$\cos\th$, so it suffices to show that the right hand side is continuous as
well. For this we have to show that the convergence of the right hand side
is uniform with respect to $\cos\th$. This follows from the estimates
$$
\gather
\Bigl\vert {{(q^{1+n};q)_\infty}\over{(q;q)_\infty}} \,
{}_1\vp_1\left({{q^{-m}}\atop{q^{1+n}}};q,a^2q^{m+n-\nu}z\right)\Bigr\vert
\leq (-\vert a^2z\vert q^{n-\n};q)_\infty \\
\Bigl\vert {{(q^{\n+n+1};q)_\infty}\over{(q;q)_\infty}} \,
{}_2\vp_1\left({{abq^{n+m},0}\atop{q^{\n+n+1}}};q,z\right)\Bigr\vert \leq
{{(-q^{\Re \n +1},-\vert ab\vert;q)_\infty}\over{(q,\vert z\vert;q)_\infty}}
\endgather
$$
and the asymptotic behaviour of the Al-Salam--Chihara polynomials given by
$$
S_n(\hf (\xi+\xi^{-1});a,b\mid q) = \xi^{-n} \, A(\xi) + {\Cal O}(\xi^{-n}),
\quad n\to\infty,\ \vert \xi\vert <1,
\tag\eqname{\vglasymASCp}
$$
with $A(\xi)= (a\xi,b\xi;q)_\infty/(\xi^2;q)_\infty$ off the spectrum and by
$$
\gather
S_n(\cos\th;a,b\mid q) = 2\vert A(e^{i\th})\vert \cos(n\th-\phi) +
{\Cal O}(q^{n/2}),\quad n\to\infty,\ 0<\th<\pi,\ \phi=\arg A(e^{i\th}),\\
S_n(\pm 1;a,b\mid q) = (\pm 1)^n n
{{(\pm a,\pm b;q)_\infty}\over{(q;q)_\infty}} + {\Cal O}(1), \quad n\to\infty
\endgather
$$
on the spectrum, cf. Askey and Ismail \cite{\AskeI, \S 3.1}, Ismail and
Wilson \cite{\IsmaW, (1.11), (1.13), \S 3}.


\subhead\newsection . Special and limiting cases\endsubhead
The first special case of interest of the addition and product formula is the
case $m=0$, which gives the decomposition of the ${}_2\vp_1$-series involved
in terms of Al-Salam--Chihara polynomials. Some other special and limiting
cases are described in the rest of this section.

\subsubhead The Koornwinder-Swarttouw addition formula as a limit
case\endsubsubhead
Formally we can obtain the $q$-analogue of Graf's addition formula
for the Jackson $q$-Bessel function derived by Koornwinder and Swarttouw
\cite{\KoorS, (4.10)} as a
special case of the addition formula \thetag{\vgladditionform} by use of the
asymptotic behaviour of the Al-Salam--Chihara polynomials off the spectrum,
cf. \thetag{\vglasymASCp}, as follows. Let $m\to\infty$ in
\thetag{\vgladditionform} and use that formally
$$
\lim_{m\to\infty} {{S_{n+m}(\hf(\xi+\xi^{-1});a,b\mid
q)}\over{S_m(\hf(\xi+\xi^{-1});aq^{-\n},b\mid q)}} = \xi^{-n}
{{(a\xi;q)_\infty}\over{(aq^{-\n}\xi;q)_\infty}},\quad \vert \xi\vert < 1,
$$
by \thetag{\vglasymASCp}. Replace in the resulting formula $z$, $a$, and
$\xi$ by $-y^2$, $q^{\n/2}x/y$, and $q^{\n/2}/s$ to obtain formally
$$
\multline
y^\n {{(xy^{-1}s^{-1},q^{\n+1};q)_\infty}\over{
(q^\n xy^{-1}s^{-1},q;q)_\infty}} \, {}_2\vp_1 \left( {{q^\n xy^{-1}s^{-1},
xy^{-1}s}\atop{q^{\n+1}}};q,-y^2\right) =
\sum_{n=-\infty}^\infty s^n y^{\n+n}  \\ \times
{{(q^{\n+n+1};q)_\infty}\over{(q;q)_\infty}} \, {}_2\vp_1\left(
{{0,0}\atop{q^{\n+n+1}}};q,-y^2\right)
x^n q^{\hf n(n-1)}
{{(q^{n+1};q)_\infty}\over{(q;q)_\infty}} \, {}_0\vp_1\left(
{{-}\atop{q^{n+1}}};q,-x^2q^n\right),
\endmultline
\tag\eqname{\vglKoorSaddf}
$$
which has been proved rigorously by Koornwinder and Swarttouw \cite{\KoorS,
(4.10)} using generating function techniques for $\n\in\Z$.

\subsubhead The limit $q\uparrow 1$\endsubsubhead
In \thetag{\vglKoorSaddf} we replace $x$ and $y$ by $(1-q)x$ and $(1-q)y$
before we take the limit $q\uparrow 1$. If we use the $q$-gamma function
$\Gamma_q(x) = (q;q)_\infty(q^x;q)_\infty^{-1}(1-q)^{1-x}$,
cf. \cite{\GaspR, \S 1.10}, and the limit
relation $\Gamma_q(x)\to \Gamma(x)$ as $q\uparrow 1$, we see that
\thetag{\vglKoorSaddf} goes over into Graf's addition
formula \thetag{\vglGrafadditie}.
Since \thetag{\vglKoorSaddf} is a limiting case of
\thetag{\vgladditionform}, we have shown that \thetag{\vgladditionform}
is a $q$-analogue of Graf's addition formula \thetag{\vglGrafadditie}.

It is also possible to use the techniques of Van Assche and Koornwinder
\cite{\VanAK, thm.~1} to treat the limit case $q\uparrow 1$ of the addition
formula \thetag{\vgladditionform} to Graf's addition formula
\thetag{\vglGrafadditie}. To this end observe that from
\thetag{\vglccfASCpols} we have
$$
{{S_{n+m}(x;a,b\mid q)}\over{S_m(x;aq^{-\n},b\mid q)}} =
\sum_{k=0}^{n+m} {{(q^\n;q)_k}\over{(q;q)_k}} (aq^{-\n})^k
(q^{m+n-k+1};q)_k
{{S_{n+m-k}(x;aq^{-\n},b\mid q)}\over{S_m(x;aq^{-\n},b\mid q)}}.
$$
In the right hand side we replace $q=c^{1/m}$, $c\in (0,1)$, and let
$m\to\infty$. Then \cite{\VanAK, thm.~1} can be used to evaluate this limit
and the rest of the limit transition is a straightforward exercise using the
binomial formula. See also
\cite{\KoelS, \S 4} for the details of a similar limit transition.

The limit transition of the product formula \thetag{\vglproductform} to
the product formula \thetag{\vglGrafproduct} for the Bessel function as
$q\uparrow 1$ is treated by use of theorem~2 of Van Assche en Koornwinder
\cite{\VanAK}. This time we have to use the connection coefficient
formula \thetag{\vglccfASCpols} in the form
$$
S_m(x;aq^{-\n},b\mid q) = \sum_{k=0}^m {{(q^{-\n};q)_k}\over{(q;q)_k}}
a^k (q^{m-k+1};q)_k S_{m-k}(x;a,b\mid q)
$$
before replacing $q=c^{1/m}$, $c\in (0,1)$, and letting $m\to\infty$.
Invoking \cite{\VanAK, thm.~2} and the binomial theorem shows that
\thetag{\vglproductform} tends \thetag{\vglGrafproduct}. See also
\cite{\KoelS, \S 5} for the details of a similar limit transition.

\subsubhead Orthogonality relations for $q$-Charlier
polynomials\endsubsubhead
The addition formula \thetag{\vgladditionform}
is a $q$-analogue of Graf's addition formula for the Bessel function
from which the Hansen-Lommel orthogonality relations for the Bessel
functions, $\sum_{k=-\infty}^\infty J_n(z)J_{n+p}(z)=\d_{0,p}$, $p\in\Z$,
$z\in\C$, can be derived, cf. Watson
\cite{\Wats, \S\S 2.4, 2.5, 11.2, 11.3}. Here we can also
specialise the parameters to obtain orthogonality relations from the addition
formula. Take $a=b=e^{i\th}=q^\hf$, $\n=p\in\Z$ and observe that
$$
\gather
{{(q^{p+1};q)_\infty}\over{(q;q)_\infty}}\, {}_2\vp_1 \left(
{{1,q}\atop{q^{p+1}}};q,z\right) = \cases 0, &\text{$p<0$,}\\
(q;q)_p^{-1}, &\text{$p\geq 0$,}\endcases  \\
S_m(\hf(q^\hf+q^{-\hf});q^{\hf-p},q^\hf\mid q) =q^{-\hf m} (q^{1-p};q)_m
= \cases (q;q)_mq^{-\hf m}, &\text{$p=0$,} \\
0, &\text{$0<p\leq m$,}\\
(q^{1-p};q)_m q^{-\hf m}, &\text{$p>m$,}\endcases \\
S_{n+m}(\hf(q^\hf+q^{-\hf});q^\hf,q^\hf\mid q) =q^{-\hf (n+m)}
(q;q)_{n+m}
\endgather
$$
to find for $p\leq m$ the orthogonality relations
$$
\multline
\d_{0,p} (q;q)_m =\sum_{n=-m}^\infty (-z)^n q^{\hf n(n-1)}(q;q)_{n+m}
{{(q^{n+1};q)_\infty}\over{(q;q)_\infty}}\, {}_1\vp_1 \left(
{{q^{-m}}\atop{q^{n+1}}};q,q^{1+m+n-p}z\right) \\ \times
{{(q^{n+p+1};q)_\infty}\over{(q;q)_\infty}}\, {}_2\vp_1 \left(
{{q^{n+m+1},0}\atop{q^{n+p+1}}};q,z\right).
\endmultline
\tag\eqname{\vglHLorthrel}
$$
Repace $z$ by $-z^2q^p/4$ and let $m\to\infty$ in \thetag{\vglHLorthrel}
to find Hansen-Lommel orthogonality relations for the Jackson $q$-Bessel
function, cf. \cite{\KoelJMAA, thm.~3.1}, which can also be obtained
from the Koornwinder-Swarttouw $q$-analogue of Graf's addition formula
\thetag{\vglKoorSaddf}, cf. \cite{\KoelJMAA, rem.~1, p.~432}.

In the orthogonality relations \thetag{\vglHLorthrel} we
use the limiting case $b\to 0$ of Heine's
transformation formula \cite{\GaspR, (1.4.6)}, cf.
$$
{}_2\vp_1\left( {{a,0}\atop{c}};q,z\right) = {1\over{(z;q)_\infty}} \,
{}_1\vp_1\left( {{c/a}\atop{c}};q,az\right).
\tag\eqname{\vgltransf}
$$
Next we replace $n$, $p$ and $z$ by $h-m$, $m-r$, $-aq^{-r}$ to see that
\thetag{\vglHLorthrel} is equivalent to the orthogonality relations for the
$q$-Charlier polynomials, cf. \cite{\GaspR, ex.~7.13, with the squared norm
replaced by its reciprocal},
$$
\sum_{h=0}^\infty {{a^h q^{\hf h(h-1)}}\over{(q;q)_h}} \bigl( c_mc_r\bigr)
(q^{-h};a;q) = \d_{m,r} q^{-m} (-a^{-1}q, q;q)_m (-a;q)_\infty,
\tag\eqname{\vglorthoqCharlier}
$$
where the $q$-Charlier polynomials are defined by
$$
\align
c_m(x;a;q) &= \, {}_2\vp_1\left( {{q^{-m},x}\atop{0}};q,-{{q^{m+1}}\over
a}\right) \\ &= (-a)^{-m} q^{m^2} x^m (x^{-1}q^{1-m};q)_m \,
{}_1\vp_1\left( {{q^{-m}}\atop{x^{-1}q^{1-m}}};q,-{{aq^{1-m}}\over x}\right).
\endalign
$$
The last equality follows by series inversion. So the $q$-Charlier polynomials
are related to Moak's $q$-Laguerre polynomials \cite{\Moak} by
$$
\aligned
c_m(q^{-\a-m};a;q) &= (-aq^\a)^{-m} (q;q)_m L_m^{(\a)} (aq^m;q), \\
L_n^{(\a)}(x;q) &= {{(q^{\a+1};q)_n}\over{(q;q)_n}}\, {}_1\vp_1\left(
{{q^{-n}}\atop{q^{\a+1}}};q,-xq^{\a+n+1}\right).
\endaligned
\tag\eqname{\vglqChaqLag}
$$

The orthogonality relations
\thetag{\vglorthoqCharlier} for the $q$-Charlier polynomials have been
obtained from the quantum algebra approach by Kalnins, Miller and Mukherjee
\cite{\KalnMMeen, (3.2)} and Floreanini and Vinet \cite{\FlorVtwee, (59)}
using representations of the $q$-oscillator algebra.
The orthogonality relations are also a byproduct of the quantum group
theoretic proof of an addition formula for the big $q$-Legendre polynomial,
cf. \cite{\KoelCJM, cor.~4.2}. In the limit case $m\to\infty$ of
\thetag{\vglHLorthrel} we know that the dual orthogonality relations also
hold, cf. \cite{\KoelJMAA, thm.~3.3}, but the orthogonality relations dual
to the orthogonality relations \thetag{\vglorthoqCharlier} for the
$q$-Charlier polynomials do not hold, cf. \cite{\KoelCJM, prop.~4.1,
cor.~4.2}. This is not correct in \cite{\FlorVtwee, (60)}.

\subsubhead More identities for $q$-Charlier polynomials\endsubsubhead
The transformation formula \thetag{\vgltransf} can also be applied in the
general addition formula \thetag{\vgladditionform}. If we next replace
$a$, $b$, $z$, $n$, and $\n$ by $q^{\hf(\m+1)}\sqrt{\a/\b}$,
$q^{\hf(\m+1)}\sqrt{\b/\a}$, $-\b q^{-r}$, $h-m$, and $m-r+\m$, we get the
following extension of the orthogonality relations
for the $q$-Charlier polynomials;
$$
\multline
\sum_{h=0}^\infty {{(\a\b)^{\hf h} q^{\hf h(\m+1)} q^{\hf
h(h-1)}}\over{(q;q)_h}} {{(q^{1+h+\m};q)_\infty}\over{(q;q)_\infty}}
c_m(q^{-h};\a;q) c_r(q^{-h-\m};\b;q) \\ \times S_h(\cos\th;q^{\hf(\m+1)}
\sqrt{ \a/\b};q^{\hf(\m+1)}\sqrt{\b/\a}\mid q) =
(-1)^{m+r} q^{\hf m(m+\m)} q^{r(r-m-\m)} \a^{-\hf m} \b^{\hf m -r} \\
\times(-\b q^{-r};q)_\infty {{(q^{1+m-r+\m};q)_\infty}\over{(q;q)_\infty}}
\, {}_2\vp_1\left( {{q^{\hf (\m+1)}\sqrt{\a/\b}e^{i\th},
q^{\hf (\m+1)} \sqrt{\a/\b}e^{-i\th}}\atop{q^{1+m-r+\m}}};q,-\b q^{-r}\right)
\\ \times S_m(\cos\th;q^{r-m+ \hf(1-\m)}\sqrt{\a/\b};
q^{\hf(\m+1)}\sqrt{ \b/\a}\mid q).
\endmultline
$$
Using \thetag{\vglqChaqLag} we can also rewrite this as an identity for
$q$-Laguerre polynomials, which then gives an alternative for the
addition formulas for $q$-Laguerre polynomials derived by Kalnins, Manocha
and Miller \cite{\KalnMaMi, (7.14)} and Kalnins and Miller
\cite{\KalnM, (4.13)} using representations of the $q$-oscillator algebra.

For particular choices of $\cos\th$ it is possible to evaluate the
Al-Salam--Chihara polynomials, which then simplifies the formula for
the $q$-Charlier polynomials. In particular, for
$e^{i\th}=q^{\hf (\m+1)}\sqrt{\b/\a}$ the Al-Salam--Chihara
polynomials reduce to $q$-shifted factorials as in the proof of
\thetag{\vglHLorthrel} and we obtain
$$
\multline
\sum_{h=0}^\infty {{\a^h q^{\hf h(h-1)}}\over{(q;q)_h}}
c_m(q^{-h};\a;q) c_r(q^{-h-\m};\b;q) =
(-1)^{m+r} q^{\hf m(m-1)} q^{r(r-m-\m)} \b^{-r} \\
\times (q^{1+r-m};q)_m (-\b q^{-r};q)_\infty
{{(q^{1+m-r+\m};q)_\infty}\over{(q^{1+\m};q)_\infty}}
\, {}_2\vp_1\left( {{q^{\m+1},\a/\b}\atop{q^{1+m-r+\m}}};q,-\b q^{-r}\right)
\endmultline
$$
as a nice extension of the orthogonality relations of the $q$-Charlier
polynomials.


\subhead\newsection . Concluding remark\endsubhead
As already said in the introduction, the addition formula
\thetag{\vgladditionform} is originally derived from the interpretation of
certain $q$-Bessel functions on the quantum group of plane motions.
Actually, the addition formula follows from
the identity \cite{\KoelIM, (4.6) with $r=\infty$}, which reflects the
homomorphism property of a representation of a group which classically gives
addition formulas if the matrix elements are known in terms of special
functions. As indicated in
remark~4.3 of \cite{\KoelIM} the elements in this identity
are explicitly (but only on a formal
level) known in terms of a certain non-commutative algebra $\A$.
Then we apply the infinite dimensional $\tau\otimes\tau$ as defined in
\cite{\KoelDMJ, \S 5} to \cite{\KoelIM, (4.6) with $r=\infty$} to obtain an
operator identity in $\ell^2(\Z)\otimes\ell^2(\Z)$. By letting it act on
suitable vectors in the representation space we can convert
the operator identity to
the addition formula \thetag{\vgladditionform} for special values of
$a$, $b$ and $\n$.

A certain $q$-analogue of Graf's addition formula for the Hahn-Exton
$q$-Bessel function can also be proved from a similar quantum group
theoretic interpretation, cf. \cite{\KoelDMJ, \S 6}. (See
\cite{\KalnMMtwee, \S 3} for a quantum algebra theoretic proof and
\cite{\KoelS} for an analytic proof as well as for
the limit case $q\uparrow 1$.)
In this case the starting point is \cite{\KoelDMJ, (6.1)}, which can be
considered as the limit case $s,t\to\infty$ of the starting point
\cite{\KoelIM, (4.6) with $r=\infty$} for the addition formula in this
paper. So it is tempting to think that the $q$-analogue
of Graf's addition
formula for the Hahn-Exton $q$-Bessel functions (\cite{\KoelDMJ,
thm.~6.3}, \cite{\KalnMMtwee, \S 3}, \cite{\KoelS, (1.4)})
can be obtained as a
limiting case of the addition formula \thetag{\vgladditionform}.
However, I have not been able to show this.


\enddocument